\def\b#1{{\bf #1}}
\def\i#1{{\it #1}}

\newtheorem{definition}{Definition}

\documentstyle[11pt]{article}
\begin{document}

\title{Polymatrix and Generalized Polynacci Numbers}
\author{Mario Catalani\\
Department of Economics, University of Torino\\ Via Po 53, 10124 Torino, Italy\\
mario.catalani@unito.it}
\date{}
\maketitle
\begin{abstract}
\small{We consider $m$-th order linear recurrences that can be thought of
as generalizations of the Lucas sequence. We exploit some interplay with
matrices that again can be considered generalizations of the Fibonacci
matrix. We introduce the definition of reflected sequence and inverted
sequence and we establish some relationship between the coefficients
of the Cayley-Hamilton equation for these matrices and the introduced
sequences.}
\end{abstract}

\section{Antefacts}
Let us define the $m\times m$ matrix $\b{A}_m$ as a matrix with the first
column of all ones, as well as the first upper diagonal, while all the other
elements are equal to zero, that is
$$\left \{\begin{array}{lccl}a_{i1}^{(m)}&=&1&\mbox{for $i=1,\ldots,m$}\\
a_{i,i+1}^{(m)}&=&1&\mbox{for $i=1,\ldots,m-1$}\\
a_{ij}^{(m)}&=&0&\mbox{otherwise.}\end{array}\right .$$
We can write in partitioned form
\begin{equation}
\label{eq:partizione}
\b{A}_m=\left[\begin{array}{cc}\b{A}_{m-1}&\b{e}_1\\
\b{e}_2'&0\end{array}\right ],
\end{equation}
where $\b{e}_1$ and $\b{e}_2$ are $m\times 1$ vectors
$$\b{e}_1=\left [\begin{array}{c}0\\0\\ \vdots\\ 0\\1\end{array}\right ],
\qquad
\b{e}_2\left [\begin{array}{c}1\\0\\ \vdots\\ 0\\0\end{array}\right ].$$

\noindent
Note that the inverse of $\b{A}_m$ is a matrix $\b{B}_m$ such that
$$\left \{\begin{array}{lccl}b_{1m}&=&1\\
b_{im}&=&-1&\mbox{for $i=2,\ldots,m$}\\
b_{i,i-1}&=&1&\mbox{for $i=2,\ldots,m$}\\
b_{ij}&=&0&\mbox{otherwise.}\end{array}\right .$$
If
we evaluate the determinant developping according the elements of the
last row we see that the only non-zero summand corresponds to element in
position $(m,1)$ while the corresponding minor is a $m-1\times m-1$
identity matrix: it follows
$$\vert\b{A}_m\vert =(-1)^{m+1}.$$

\noindent
The characteristic polynomial of $\b{A}_m$ is
\begin{equation}
\label{eq:caratteristico}
g(x)=x^m-x^{m-1}-x^{m-2}-\cdots -x-1.
\end{equation}
Now, just to fix notation, given $n$ numbers
$$\{a_i\}_1^n=\{a_1,\,a_2,\,\ldots ,\,a_n\}$$
let $S_k^n\left <a_i\right >$ denote the symmetric functions, that is
$$S_k^n\left <a_i\right >=
\sum_{\stackrel{1\le i_1,\ldots ,i_k\le n}{i_1<i_2<\cdots <i_k}}
a_{i_1}a_{i_2}\cdots a_{i_k}, \qquad 1\le k\le n.$$
The Cayley-Hamilton equation for $\b{A}_m$ is
$$\b{A}_m^m-\b{A}_m^{m-1}-\cdots -\b{A}_m-\b{I}=\b{0}.$$
Using relationships among coefficients of the Cayley-Hamilton equation and
the eigenvalues $\{r_i,\,i=1,\ldots ,m\}$ we have
$$S_i^m(r_i)=(-1)^{i+1},\qquad i=1,\ldots ,m.$$
Note that the maximal real root approaches
2 for $m$ going to infinity (see \cite{wolfram}).

\noindent
From the characteristic polynomial we can
define the $m$-th order recurrence
\begin{equation}
\label{eq:ricorrenza1}
U_n^{(m)}=U_{n-1}^{(m)}+U_{n-2}^{(m)}+\cdots +U_{n-m}^{(m)},
\quad n\ge m.
\end{equation}
With $m=2$ we have the Fibonacci-Lucas sequence, with $m=3$ the Tribonacci
sequence, with $m=4$ the Tetranacci sequence, and so on. For this
reason we might call this matrix a Polynacci matrix (Polymatrix).

\noindent
Note that (see \cite{mario4})
$$\sum_{i_1=0}^n \sum_{i_2=0}^{n-i_1}\cdots \sum_{i_{m-1}=0}^{n-i_1-\cdots
-i_{m-2}}r_1^{i_1}r_2^{i_2}\times\cdots\times
r_{m-1}^{i_{m-1}}r_m^{n-i_1-\cdots -i_{m-1}}
=V_{n+1},$$
where
$$V_n=V_{n-1}+V_{n-2}+\cdots +V_{n-m+1}+V_{n-m},$$
with initial conditions
$$V_0=\,V_1=V_2=\cdots =V_{m-2}=V_{m-1}=1.$$

\section{Generalized Polynacci Sequences}
Now we want to choose the initial conditions of $U_n^{(m)}$ in such a way
that, for any $m$, there holds the Binet form
$$U_n^{(m)}=r_1^n+r_2^n+\cdots +r_m^n.$$
In this way we obtain generalized Polynacci sequences. The term
\i{generalized} stems from the fact that the Tribonacci numbers ($m=3$) so
defined bears with the Tribonacci with initial conditions $0,\,1,\,1$
(see \cite{koshy}) the same resemblance as the Lucas sequence does with
the Fibonacci sequence.

\noindent
It follows
$$U_n^{(m)}={\rm tr}(\b{A}_m^n).$$
Now we are going to establish the following result
\begin{equation}
\label{eq:fondamentale}
{\rm tr}(\b{A}_m^k)={\rm tr}(\b{A}_{m-1}^k), \qquad k<m,
\end{equation}
where ${\rm tr}(\cdot)$ is the trace operator,
so that
\begin{equation}
{\rm tr}(\b{A}_m^k)={\rm tr}(\b{A}_k^k).
\end{equation}
We assume $k>0$, since ${\rm tr}(\b{A}_m^0)=m$.
We need the following result, which can be established with tedious
calculations,
\begin{equation}
\label{eq:scalare}
\b{e}_2'\b{A}_m^i\b{e}_1=0, \qquad i=0,\,1,\,\ldots,\,m-2.
\end{equation}
Using repeated matrix multiplication on the partitioned form
(Equation~\ref{eq:partizione}) of $\b{A}_m$ the asked for trace is the
sum of the trace of the resulting first diagonal blok plus the resulting
scalar which is the second diagonal block. Starting analyzing this scalar
we see that, using Equation~\ref{eq:scalare}, for $h\le k$, this term in
$\b{A}_m^h$ is $\b{e}_2'\b{A}_{m-1}^{h-2}\b{e}_1$. Since $h\le k$ and
$k<m$ it follows $h-2<m-2$ so that using again Equation~\ref{eq:scalare}
with $m$ replaced by $m-1$ we have
$\b{e}_2'\b{A}_{m-1}^{h-2}\b{e}_1=0,\;\forall\, h.$

\noindent
As for the first diagonal block the resulting expression in $\b{A}_m^k$
is a rather messy sum one term of which is $\b{A}_{m-1}^k$. But since we
are interested just in the trace, using the fact that for a matrix \b{B}
and vectors \b{u} and \b{z} we have ${\rm tr}(\b{Auz}')=
{\rm tr}(\b{z}'\b{Au})$, we
get that the trace of this block is ${\rm tr}(\b{A}_{m-1}^k)$ plus the sum of
terms of the form
\begin{equation}
\label{eq:casino}
c_i\prod_{j=1}^{\left\lceil{k\over 2}\right\rceil-1}
\b{e}_2'\b{A}_{m-1}^{2j+\alpha}\b{e}_1,
\end{equation}
where
$$\alpha=2\left\{\left\lfloor
{k\over 2}\right\rfloor-{k\over 2}\right\},$$
that is $\alpha=0$ if $k$ is even,
and $\alpha=-1$ if $k$ is odd, and $c_i$ are constant. Since the greatest
power of $\b{A}_{m-1}$ in Equation~\ref{eq:casino} is $k-2$ and $k-2<m-2$
again invoking Equation~\ref{eq:scalare} all these terms turn out to be
equal to zero. From this we obtain
Equation~\ref{eq:fondamentale}.

\noindent
It follows that
we can determine the initial conditions recursively
from Lucas initial conditions ($m=2$), obtaining,
for example, for generalized Tribonacci
(A001644 in \cite{sloane})
3,1,3; for generalized Tetranacci (A073817 in \cite{sloane})
4,1,3,7; for generalized Pentanacci (A074048 in \cite{sloane}) 5,1,3,7,15; for
generalized Hexanacci (A074584 in \cite{sloane}) 6,1,3,7,15,31.

\noindent
But we can also obtain
a general formula. As we know $U_0^{(m)}=m,\,U_1^{(m)}=1$. We are going to
prove that
\begin{equation}
\label{eq:dimostrazione}
U_i^{(m)}=2^i-1,\quad i=2,\,\ldots,\,m-1.
\end{equation}
This will be done easily by induction on $m$. This is true for $m=2$ since
$U_2^{(2)}$ is $L_2$ which is $3=2^2-1$. Assume that the claim holds for
$m$: then we have to show that
$$U_i^{(m+1)}=2^i-1,\quad i=2,\,\ldots,\,m.$$
Now, for $i=2,\,\ldots,\,m-1$,
\begin{eqnarray*}
U_i^{(m+1)}&=&{\rm tr}(\b{A}_{m+1}^i)\\
&=&{\rm tr}(\b{A}_m^i)\\
&=& U_i^(m)\\
&=&2^i-1,
\end{eqnarray*}
because of Equation~\ref{eq:fondamentale}. For $i=m$ we have
\begin{eqnarray*}
U_m^{(m+1)}
&=&{\rm tr}(\b{A}_{m+1}^m)\\
&=&{\rm tr}(\b{A}_m^m)\\
&=&U_m^{(m)}.
\end{eqnarray*}
But $U_m^{(m)}$ is the first determined by the recurrence, so
\begin{eqnarray*}
U_m^{(m)} &=& \sum_{i=0}^{m-1}U_i^{(m)}\\
&=&m+1+\sum_{i=2}^{m-1}(2^i-1)\\
&=&2^m-1,
\end{eqnarray*}
since
$$\sum_{i=2}^{m-1}(2^i-1)=2^m-m-2,$$
and so the claim is proved.

\noindent
Now we are going to derive a closed form of the ordinary generating
function (ogf). So let $G(x)$ be the ogf of $U_n^{(m)}$
$$G(x)= a_0+a_1x+a_2x^2+a_3x^3+ \cdots $$
where $a_i=U_i^{(m)}$.
Multiply both sides by $x^i,\,i=1,\,\ldots,\,m$. Then
\begin{eqnarray*}
&&G(x)(1-x-x^2-x^3-\ldots -x^m)=a_0+(a_1-a_0)x+(a_2-a_1-a_0)x^2\\
&&\quad +(a_3-a_2-a_1-a_0)x^3+\cdots +(a_m-a_{m-1}-\cdots -a_1-a_0)x^m
+\cdots
\end{eqnarray*}
Because of the recurrence relationship all the coefficients of $x^j,\,
j\ge m$ are equal to zero. Now insert the initial conditions $a_i,\,
i=0,\ldots, m-1$ and we are left with
\begin{eqnarray*}
&&G(x)(1-x-x^2-x^3-\ldots -x^m)=m+(1-m)x+(2-m)x^2+\cdots +\\
&&\qquad\qquad\qquad\qquad
+\left (2^{m-1}-1-\sum_{i=2}^{m-2}(2^i-1)-1-m\right )x^{m-1}.
\end{eqnarray*}
The last summand turns out equal to $-1$. Then we can conclude
\begin{equation}
G(x)={m-(m-1)x-(m-2)x^2-\cdots -x^{m-1}\over
1-x-x^2-x^3-\cdots -x^m}.
\end{equation}
Note that we can allow for negative subscripts, following
and generalizing \cite{howard}. If we have a general $m$-th order
linear recurrence
\begin{equation}
U_n^{(m)}=s_1U_{n-1}^{(m)}+s_2U_{n-2}^{(m)}+\cdots +s_mU_{n-m}^{(m)},
\quad n\ge m
\end{equation}
then we can define
\begin{equation}
\label{eq:negativo}
U_{-n}^{(m)}=-{s_{m-1}\over s_m}U_{-(n-1)}^{(m)}-{s_{m-2}\over s_m}
U_{-(n-2)}^{(m)}-\cdots -{s_1\over s_m}U_{-(n-m+1)}^{(m)}
+{1\over s_m}U_{-(n-m)}^{(m)}.
\end{equation}

\section{Reflected Sequences}
The reflected polynomial (\cite[p. 339]{knuth}) of $g(x)$ is
\begin{equation}
\label{eq:reflected}
g^R(x)=-x^m-x^{m-1}-x^{m-2}-\cdots -x+1.
\end{equation}
The roots are the reciprocals of the roots of $g(x)$, that is $\left\{
{1\over r_i}\right\}$, and so $g^R(x)$ is the characteristic polynomial
of matrix $\b{B}_m$.

\begin{definition} The \b{reflected} recurrence ${\tilde U}_n^{(m)}$
of recurrence~\ref{eq:ricorrenza1}
is the recurrence with characteristic polynomial which is the reflected
characteristic polynomial and with initial conditions such that the
coefficients of the respective Binet forms are the same.
\end{definition}

It follows
$${\tilde U}_n^{(m)}=-{\tilde U}_{n-1}^{(m)}-{\tilde U}_{n-2}^{(m)}
-\cdots -{\tilde U}_{n-m+1}^{(m)}+{\tilde U}_{n-m}^{(m)},\quad n\ge m,$$
and
$${\tilde U}_n^{(m)}={1\over r_1^n}+{1\over r_2^n}+\cdots
+{1\over r_m^n}={\rm tr}(\b{B}_m^n).$$
Now we are going to evaluate ${\rm tr}(\b{B}_m^n)$, for $n=0,\,1,\,
\ldots\, m-1$ so that we get the required initial conditions. Of
course ${\rm tr}(\b{B}_m^0)=m$. We know what is the expression for
$\b{B}_m$, so ${\rm tr}(\b{B}_m^1)=-1$. Now consider what happens when
we perform matrix multiplication, starting from $\b{B}_m\b{B}_m$: when we
postmultiply by $\b{B}_m$ then the columns of the matrix to the left are
shifted to the left by one place and the last column is the linear
combination of the columns of the matrix to the left with coefficients
$\{1,-1,-1,\ldots ,-1\}$. It is easy to see that for $\b{B}_m\b{B}_m$ the
last column is $\{-1,2,0,\ldots ,0\}$. Then all diagonal elements are
zero except for the element in the next-to-last column which is $-1$, so
that the trace is $-1$. If we go ahead and consider $\b{B}_m^3$ we see that
the last column is $\{0,-1,2,0,\ldots ,0\}$ and the only non zero diagonal
element appears in column $m-2$ and it is equal to $-1$, so the trace is $-1$.
Repeating the same reasoning we come to $\b{B}_m^{m-1}$: here the last
column is $\{0,\ldots, 0,-1,2,0\}$ and the only non zero diagonal element
appear in column 2 and it is equal to $-1$, so the trace is $-1$. So
we have proved that the asked for initial conditions are
$${\tilde U}_0^{(m)}=m,\, {\tilde U}_1^{(m)}={\tilde U}_2^{(m)}=\cdots
={\tilde U}_{m-1}^{(m)}=-1.$$
Of course either using the recurrence or going on with the multiplication
process it turns out that ${\tilde U}_m^{(m)}=2m-1$.
If for example $m=4$ we get
$${\tilde U}_0^{(4)}=4,\,{\tilde U}_1^{(4)}=-1,\,{\tilde U}_2^{(4)}=-1,\,
{\tilde U}_3^{(4)}=-1,$$
while
$$U_0^{(4)}=4,\, U_1^{(4)}=1,\, U_2^{(4)}=3,\,
U_3^{(4)}=7.$$
Now we are going to derive a closed form of the ordinary generating
function (ogf) for the reflected recurrence.
So let ${\tilde G}(x)$ be the ogf of ${\tilde U}n^{(m)}$
$${\tilde G}(x)= a_0+a_1x+a_2x^2+a_3x^3+ \cdots $$
where $a_i={\tilde U}_i^{(m)}$.
Multiply both sides by $x^i,\,i=1,\,\ldots,\,m$. Then
\begin{eqnarray*}
&&\!\!\!\!\!\!\!\!\!\!\!{\tilde G}(x)(1+x+x^2+-\ldots +x^{m-1}-x^m)=
a_0+(a_1+a_0)x+(a_2+a_1+a_0)x^2\\
&&\quad +(a_3+a_2+a_1+a_0)x^3+\cdots +(a_m+a_{m-1}+\cdots +a_1-a_0)x^m
+\cdots
\end{eqnarray*}
Because of the recurrence relationship all the coefficients of $x^j,\,
j\ge m$ are equal to zero. Now insert the initial conditions $a_i,\,
i=0,\ldots, m-1$ and we are left with
\begin{eqnarray*}
&&{\tilde G}(x)(1+x+x^2+\ldots +x^{m-1}-x^m)=m+(m-1)x+(m-2)x^2+\cdots \\
&&\qquad\qquad\qquad\qquad\qquad
+(2m-1-1-\cdots -1-m)x^{m-1}.
\end{eqnarray*}
The last summand turns out equal to $1$. Then we can conclude
\begin{equation}
{\tilde G}(x)={m+(m-1)x+(m-2)x^2+\cdots +x^{m-1}\over
1+x+x^2+\cdots +x^{m-1}-x^m}.
\end{equation}
Note that for generalized Polynacci sequences Equation~\ref{eq:negativo}
becomes
\begin{equation}
U_{-n}^{(m)}=-U_{-(n-1)}^{(m)}-
U_{-(n-2)}^{(m)}-\cdots - U_{-(n-m+1)}^{(m)}
+U_{-(n-m)}^{(m)}.
\end{equation}
From Equation~\ref{eq:dimostrazione} we get easily
\begin{eqnarray*}
U_{-1}^{(m)}&=&-U_0^{(m)}-
U_1^{(m)}-\cdots - U_{m-2}^{(m)}
+U_{m-1}^{(m)}\\
&=&-1.
\end{eqnarray*}
In the same way
\begin{eqnarray*}
U_{-2}^{(m)}&=&-U_{-1}^{(m)}-
U_0^{(m)}-\cdots - U_{m-3}^{(m)}
+U_{m-2}^{(m)}\\
&=&-1,
\end{eqnarray*}
and in general
$$U_{-(i-1)}=-1 \qquad i=2,\,\ldots,\,m.$$
From this it follows that
\begin{equation}
{\tilde U}_n^{(m)}= U_{-n}^{(m)}.
\end{equation}
Reflected Tribonacci is A073145, reflected Tetranacci is A074058,
reflected Pentanacci is A074062 in \cite{sloane}.

\section{Inverted Sequences}
Another related sequence ${\hat U}_n^{(m)}$ is obtained in the following way.
Define its generating function ${\hat G}(x)$ as
$${\hat G}(x)={1\over x}G\left ({1\over x}\right ).$$
Then we get
$${\hat G}(x)={1+2x+3x^2+\cdots +(m-1)x^{m-2}-mx^{m-1}\over
1+x+x^2+\cdots +x^{m-1}-x^m}.$$
Note that the numerator is the derivative of the denominator and also that
the denominator is the same as in ${\tilde G}(x)$. Using the Rational
Expansion Theorem for Distinct Roots in \cite[p. 340]{knuth} we get
easily the closed form
\begin{eqnarray}
{\hat U}_n^{(m)}&=&
-{1\over r_1^{n+1}}-{1\over r_2^{n+1}}-\cdots
-{1\over r_m^{n+1}}\nonumber\\
&=&-{\tilde U}_{n+1}^{(m)}.
\end{eqnarray}
The recurrence is then
$$
{\hat U}_n^{(m)}=-{\hat U}_{n-1}^{(m)}-{\hat U}_{n-2}^{(m)}
-\cdots -{\tilde U}_{n-m+1}^{(m)}+{\hat U}_{n-m}^{(m)},$$
with
$${\hat U}_0^{(m)}=1,\,{\hat U}_1^{(n)}=1,\,\ldots ,\,
{\hat U}_{m-1}^{(m)}=1-2m.$$
${\hat U}_n^{(3)}$ is sequence A075298 in \cite{sloane}.

\section{Cayley-Hamilton Equation}
Now consider the Cayley-Hamilton equation for $\b{A}_m^n$
$$\b{A}_m^{nm}-c_1^{(n)}\b{A}_m^{n(m-1)} +\cdots
+(-1)^ic_i^{(n)}\b{A}_m^{n(m-i)}+\cdots + (-1)^mc_m^{(n)}\b{I}=\b{0},$$
where $c_i^{(n)}$ is the sum of the determinants of the principal minors of
order $i$ of $\b{A}_m^n$. Immediately we have
$$c_1^{(n)}={\rm tr}(\b{A}_m^n)=U_n^{(m)},$$
$$c_m^{(n)}=\vert\b{A}_m\vert^n=(-1)^{(m+1)n}.$$
Note that if $m$ is odd $c_m^{(n)}=1$, while if $m$ is even $c_m^{(n)}=(-1)^n$.

\noindent
Then
\begin{eqnarray*}
{1\over r_1^n}+{1\over r_2^n}+\cdots
+{1\over r_m^n}&=&{S_{m-1}^m\left <r_i^n\right >\over \prod_{i=1}^m r_i^n}\\
&=&{S_{m-1}^m\left <r_i^n\right >\over (-1)^{(m+1)n}}\\
&=&{\rm tr}(\b{B}_m^n)\\
&=&{\tilde U}_n^{(m)}.
\end{eqnarray*}
From this it follows
\begin{equation}
\label{eq:minore}
c_{m-1}^{(n)}=(-1)^{(m+1)n}{\tilde U}_n^{(m)},
\end{equation}
that is
$c_{m-1}^{(n)}={\tilde U}_n^{(m)}$ if $m$ is odd, and
$c_{m-1}^{(n)}=(-1)^n{\tilde U}_n^{(m)}$ if $m$ is even.
Note that if the generating
function of ${\tilde U}_n^{(m)}$ is ${\tilde G}(x)$ then the
generating function
of $(-1)^n{\tilde U}_n^{(m)}$ is ${\tilde G}(-x)$.


\begin{thebibliography}{9}

\bibitem{mario4} M. Catalani (2002), "On the Roots of the Cubic Defining
the Tribonacci Sequences." arXiv:math.CO/0209265
\bibitem{knuth} R. Graham, D. Knuth and O. Patashnik (1998), \i{Concrete
Mathematics}, Addison-Wesley, Reading, MA.
\bibitem{howard}
F.T. Howard (2001), "A Tribonacci Identity." \i{The Fibonacci Quarterly}
\b{39}.4: 352-357.
\bibitem{koshy} T. Koshy (2001), \i{Fibonacci and Lucas Numbers with
Applications}, John Wiley \& Sons, New York.
\bibitem{sloane}
N.J.A. Sloane, Editor (2002),
The On-Line Encyclopedia of Integer Sequences,
http://www.research.att.com/\~{}njas/sequences/.
\bibitem{wolfram} Eric Weisstein's World of Mathematics, published
electronically at
http://mathworld.wolfram.com/Fibonaccin-StepNumber.html
\end{thebibliography}
\end{document}